\documentclass[11pt]{amsart} 
\usepackage{amssymb,amsmath,latexsym,enumerate,graphicx,bbm,mathptmx,microtype,cite}

\newtheorem{theorem}{Theorem} 

\newtheorem{exam}{Example}

\newtheorem*{rem}{Remarks}

\newcommand\commentout[1]{}
\newcommand\Def[1]{\emph{#1}}

\newcommand\ehr{\operatorname{ehr}}

\newcommand\ZZ{\mathbb{Z}}

\newcommand\RR{\mathbb{R}}

\newcommand\cO{\mathcal{O}}
\newcommand\cP{\mathcal{P}}

\newcommand\be{\mathbf{e}}

\newcommand\bx{\mathbf{x}}

\newcommand\inc{\mathrm{H}}

\begin{document}

\title{Binomial Inequalities for Chromatic, Flow, and Tension Polynomials}

\author{Matthias Beck}
\address{Department of Mathematics, San Francisco State University, U.S.A.\ \& Mathematisches
Institut, Freie Universit\"at Berlin, Germany}
\email{mattbeck@sfsu.edu}

\author{Emerson Le\'on}
\address{Departamento de Matem\'aticas \\ Universidad de los Andes\\ 
Bogot\'a\\ Colombia}
\email{emersonleon@gmail.com}


\begin{abstract}
A famous and wide-open problem, going back to at least the early 1970's, concerns the classification of chromatic polynomials of graphs.  Toward this classification problem, one may ask for necessary inequalities among the coefficients of a chromatic polynomial, and we contribute such inequalities when a chromatic polynomial
$
  \chi_G(n) = \chi^*_0 \binom {n+d} d + \chi^*_1 \binom {n+d-1} d + \dots + \chi^*_d \binom n d
$
is written in terms of a binomial-coefficient basis.  For example, we show that
$
  \chi^*_{ j } \le \chi^*_{ d-j }
$,
for $0 \le j \le \frac{ d }{ 2 }$.  Similar results hold for flow and tension polynomials
enumerating either modular or integral nowhere-zero flows/tensions of a graph.  Our theorems
follow from connections among chromatic, flow, tension, and order polynomials, as well as Ehrhart polynomials of lattice polytopes that admit unimodular triangulations.  Our results use Ehrhart inequalities due to Athanasiadis and Stapledon and are related to recent work by Hersh--Swartz and Breuer--Dall, where inequalities similar to some of ours were derived using algebraic-combinatorial methods.
\end{abstract}



\thanks{
A preliminary version of this paper \cite{leon} appeared as an extended abstract in the conference
proceedings for FPSAC 2017 (Formal Power Series and Combinatorics at Queen Mary University
of London), published in \emph{S\'eminaire Lotharingien Combinatoire}.
We thank Tristram Bogart, Tricia Hersh, Florian Kohl, Juli\'an Pulido, Alan Stapledon, Tom
Zaslavsky, and two anonymous referees
for many useful conversations and comments. Also thanks to the organizers of ECCO 2016 (Escuela
Colombiana de Combinatoria 2016 in Medell\'in, Colombia) where this research collaboration started, and to Universidad de los Andes for their support. 
}

\date{28 April 2021}

\maketitle


\section{Introduction}

A famous and wide-open problem, going back to at least~\cite{wilfchromatic}, concerns the classification of chromatic polynomials of graphs.
As is well known, for a given graph $G$, the number $\chi_G(n)$ of proper colorings of $G$ using $n$ colors evaluates to a polynomial in $n$, and so a natural question is: which polynomials are chromatic?

Toward this classification problem, one may ask for necessary inequalities among the
coefficients 
of a chromatic polynomial, and this paper gives one such set of inequalities.
In enumerative combinatorics, there are three natural bases for the space of polynomials of degree at most $d$:
\begin{itemize}
  \item the monomials $1, n, n^2, \dots, n^d$;
  \item the binomial coefficients $\binom n d, \binom n {d-1}, \dots, \binom n 0$;
  \item the binomial coefficients $\binom {n+d} d, \binom {n+d-1} d, \dots, \binom n d$.
\end{itemize}
It is well known that the coefficients  of any chromatic polynomial in the monomial basis alternate in sign (this can be proved, e.g., by
deletion--contraction), and that the coefficients in both binomial-coefficient bases are nonnegative (in the first case, this follows from
considering proper colorings that use exactly $k$ colors, for $0 \le k \le d$, and this is
closely connected to \emph{$\sigma$-polynomials}~\cite{korfhage}; in the second case, nonnegativity follows from Stanley's work
on order polynomials \cite{stanleychromaticlike} and the natural decomposition of a
chromatic polynomial into order polynomial---see equation~\eqref{eq:chiintoorder} below---; this was
first spelled out in~\cite{linialchromatic}).

We will work in the last basis and define the corresponding coefficients of the chromatic polynomial of a given graph $G$ with $d$ vertices via
\[
  \chi_G(n) \ = \ \chi^*_0 \binom {n+d} d + \chi^*_1 \binom {n+d-1} d + \dots + \chi^*_d \binom n d \, .
\]
We will collect the $\chi^*_j$s 
in the polynomial $\chi^*_G(z) := \chi^*_d \, z^d + \chi^*_{ d-1 } \, z^{ d-1 } + \dots + \chi^*_0$ (which might not
have degree $d$) and note that this polynomial appears in the generating function of $\chi_G(n)$, more precisely,
\[
  \sum_{ n \ge 1 } \chi_G(n) \, z^n \ = \ \frac{ \chi^*_G(z) }{ (1-z)^{ d+1 } } \, .
\]
To the best of our knowledge, Linial~\cite{linialchromatic} initiated the first study of the
chromatic polynomial in the form of $\chi^*_G(z)$; see
also~\cite{gansnervo,brentichromatic,tomescu,gesselchromatic}.
We think of the linear transformation going from $\chi_G(n)$
to $\chi^*_G(z)$ as a tool that is useful beyond chromatic polynomials (in fact, as we will see
below, it is a standard tool in Ehrhart theory), and so we suggest to call $\chi^*_G(z)$ the
\emph{binomial transform} of~$\chi_G(n)$.

The strongest known conditions on the coefficients\footnote{
If we do not specify a basis when talking about \emph{the coefficients} of a polynomial, we are thinking of the standard monomial basis.
} 
 of $\chi^*_G(z)$ are due to Hersh and Swartz \cite[Theorem~15]{hershswartz}: Defining $h_0, h_1, \dots, h_d$ via
\[
  \sum_{ n \ge 0 } \left( (n+1)^d - \chi_G(n+1) \right) z^n \ = \ \frac{ h_d \, z^d + h_{ d-1 } \, z^{ d-1 } + \dots + h_0 }{ (1-z)^d } \, ,
\]
they proved that
\[
  h_0 \le h_1 \le \dots \le h_{ \lfloor \frac d 2 \rfloor - 1 }
  \qquad \text{ and } \qquad
  h_j \le h_{ d-2-j } \ \text{ for } j \le \tfrac d 2 - 1 \, ,
\]
from which one can now deduce inequalities for the $\chi^*_j$s involving Eulerian numbers.

The natural dual situation concerns flows on a graph.
Denote by $\inc = (\eta_{ v,e }) \in \{ 0, \pm 1 \}^{ V \times E } $ the \emph{signed incidence matrix} of $G = (V, E)$, i.e.,
\[
  \eta_{ v,e } \ = \ \begin{cases}
    1 & \text{ if $v$ is the head of $e$, } \\
   -1 & \text{ if $v$ is the tail of $e$, } \\
    0 & \text{ if $v$ is not incident with $e$, }
  \end{cases}
\]
where we equipped $G$ with an (arbitrary but fixed) orientation.
Let $A$ be an Abelian group.
A \Def{nowhere-zero $A$-flow} on $G$ is a mapping $x: E \to A \setminus \{ 0 \}$ that is in the kernel of $\inc$.
(See, e.g., \cite{jaegernwzflowproblems,seymourhandbook} for background on nowhere-zero flows.)
Tutte \cite{tutteflowpoly} proved in 1947 that the number $\phi_G(n)$ of nowhere-zero $\ZZ_n$-flows on $G$ is a polynomial in $n$.
A more recent theorem of Kochol \cite{kocholpolynomial} says that the
number $f_G(n)$ of nowhere-zero $\ZZ$-flows on $G$ whose images satisfy $|x(e)| < n$ is also a polynomial in $n$.
(It is easy to see that both flow polynomials are independent of the chosen orientation.)
While it has long been known that $\phi_G(n)$ and $f_G(n)$ have identical integer roots, they are rather different polynomials. 

A \Def{nowhere-zero $A$-tension} on $G$ is a mapping $x: E \to A \setminus \{ 0 \}$ that is in the
row space of $\inc$.
It is not hard to see that the number of nowhere-zero $\ZZ_n$-tensions on $G$ equals $\frac{ 1 }{ n^c }
\chi_G(n)$ where $c$ denotes the number of components of $G$.
The situation for integer tensions is more interesting: Kochol \cite{kocholtension} proved that the
number $t_G(n)$ of nowhere-zero $\ZZ$-tensions on $G$ with $|x(e)| < n$ is a polynomial in $n$ (which
is quite different from~$\chi_G(n)$).

As with the chromatic polynomials, we will express $\phi_G(n)$, $f_G(n)$, and $t_G(z)$ in a binomial-coefficient basis 
(namely, $\binom {n+\xi} \xi, \binom {n+\xi-1} \xi, \dots, \binom {n} \xi$)
and define their binomial transforms 
via\footnote{
There is a subtlety here that differentiates $\chi_G(n)$ from $\phi_G(n)$, $f_G(n)$, and $t_G(z)$, and thus one needs to treat the accompanying generating functions
with some care. Namely, $\chi_G(n)$ has constant term 0, which is not true for
$\phi_G(n)$, $f_G(n)$, and $t_G(z)$. Note that we chose all of our generating functions to
start with $n=1$; the alternative choice of starting with $n=0$ would result in a
different definition of the binomial transform.
}
\[
  \sum_{ n \ge 1 } \phi_G(n) \, z^n \ = \ \frac{ \phi^*_G(z) }{ (1-z)^{ \xi + 1 } } \, ,
  \qquad 
  \sum_{ n \ge 1 } f_G(n) \, z^n \ = \ \frac{ f^*_G(z) }{ (1-z)^{ \xi + 1 } } \, ,
  \qquad \text{ and } \qquad
  \sum_{ n \ge 1 } t_G(n) \, z^n \ = \ \frac{ t^*_G(z) }{ (1-z)^{ d - c + 1 } } \, ,
\]
where $c$ denotes the number of components of $G$
and $\xi := |E| - d + c$ is the \Def{cyclomatic number} of $G$.

Similar to the chromatic situation, it is known \cite{kocholtension,kocholpolynomial,breuersanyal} that
the coefficients of $\phi^*_G (z)$, $f^*_G (z)$, and $t^*_G (z)$ are nonnegative. Breuer and
Dall~\cite{breuerdall} proved the $\ZZ_n$-flow analogues of the above-mentioned inequalities by Hersh
and Swartz: Defining $h_0, h_1, \dots, h_\xi$ via
\[
  \sum_{ n \ge 0 } \left( (n+1)^\xi - \phi_G(n+1) \right) z^n \ = \ \frac{ h_\xi \, z^\xi + h_{ \xi-1 } \, z^{ \xi-1 } + \dots + h_0 }{ (1-z)^\xi } \, ,
\]
we have
\[
  h_0 \le h_1 \le \dots \le h_{ \lfloor \frac \xi 2 \rfloor - 1 }
  \qquad \text{ and } \qquad
  h_j \le h_{ \xi-2-j } \ \text{ for } j \le \tfrac \xi 2 - 1 \, ,
\]
and again, from these one can deduce inequalities for the coefficients of $\phi^*_G(z)$ involving Eulerian numbers.
Breuer and Dall gave some constraints also for $f^*_G(z)$ and $t^*_G(z)$ but these were not as
clear cut as for~$\phi^*_G(z)$.

Our goal is to show how one can derive theorems similar to those by Hersh--Swartz and
Breuer--Dall (including inequalities for $f^*_G(n)$ and $t^*_G(z)$) through a discrete geometric setup.
Our main result is as follows.

\begin{theorem}\label{thm:chrompoldecomp}
Let $G$ be a graph on $d$ vertices with $c$ components and cyclomatic number~$\xi$.
Then 
\begin{align*}
  &\chi^*_{ 1 } \le \chi^*_{ 2 } \le \dots \le \chi^*_{ \lfloor \frac{ d+1 }{ 2 } \rfloor } \\
  &\chi^*_{ j } \le \chi^*_{ d-j } \ \text{ for } \ 1 \le j \le \tfrac {d-1} 2 \\
  &\phi^*_1 \le \phi^*_{ 2 } \le \dots \le \phi^*_{ \lfloor \frac{ \xi }{ 2 } \rfloor + 1 } \\
  &\phi^*_{ j } \le \phi^*_{ \xi+1-j } \ \text{ for } \ 1 \le j \le \tfrac \xi 2 \\
  &f^*_1 \le f^*_{ 2 } \le \dots \le f^*_{ \lfloor \frac{ \xi }{ 2 } \rfloor + 1 } \\
  &f^*_{ j } \le f^*_{ \xi+1-j } \ \text{ for } \ 1 \le j \le \tfrac \xi 2 \\
  &t^*_1 \le t^*_{ 2 } \le \dots \le t^*_{ \lfloor \frac{ d-c }{ 2 } \rfloor + 1 } \\
  &t^*_{ j } \le t^*_{ d-c+1-j } \ \text{ for } \ 1 \le j \le \tfrac {d-c} 2 \, .
\end{align*}
\end{theorem}

It is clear (though the details need some work) that the 
inequalities 
for $\chi_G^*$ and $\phi_G^*$ are closely related to the work of Hersh--Swartz and
Breuer--Dall.
At any rate, the methods in \cite{hershswartz,breuerdall} are from algebraic combinatorics: one constructs a simplicial complex whose $h$-vector satisfies certain inequalities (stemming from a convex-ear decomposition). 
Our approach, by contrast, is through Ehrhart polynomials of lattice polytopes, which we discuss in Section~\ref{sec:stapledon}.
A small subclass of lattice polytopes admit unimodular triangulations, and for this
subclass, Athanasiadis \cite{athanasiadishstareulerian} and Stapledon~\cite{stapledondelta}
proved inequalities for the binomial transforms of Ehrhart polynomials 
similar in spirit to Theorem~\ref{thm:chrompoldecomp}.

While chromatic polynomials are not Ehrhart polynomials, they can be written as sums of order polynomials (by the afore-mentioned work of Stanley~\cite{stanleychromaticlike}), which we study in
Section~\ref{sec:order}. Order polynomials, in turn, \emph{are} Ehrhart polynomials in
disguise, and so here is where we apply the Athanasiadis--Stapledon inequalities (with some
tweaking); Theorem~\ref{thm:orderpoldecomp} below might be of interest on its own accounts. 

The flow and tension inequalities in Theorem~\ref{thm:chrompoldecomp} follow in a similar fashion
from writing the (two kinds of) flow and tension polynomials as sums of Ehrhart polynomials (and then
using the Athanasiadis--Stapledon inequalities), as we illustrate in Section~\ref{sec:flow}.
For integral flows and tensions, this geometric setup was introduced by Kochol~\cite{kocholpolynomial,kocholtension}, whereas for modular flows it is due to Breuer--Sanyal~\cite{breuersanyal}.

The underlying theme here is that one can interpret graph polynomials as certain combinations
of Ehrhart polynomials of very nice polytopes (ones that admit unimodular triangulations),
and thus these Ehrhart polynomials satisfy certain linear constraints, which then can be
translated back into constraints for the graph polynomials.


\section{Ehrhart theory and $h^*$-inequalities}\label{sec:stapledon}

Given a \Def{lattice polytope} $\cP \subset \RR^d$, i.e., the convex hull of finitely many points in $\ZZ^d$, Ehrhart's celebrated
theorem~\cite{ehrhartpolynomial} says that the counting function
\[
  \ehr_\cP(n) \ := \ \left| n \, \cP \cap \ZZ^d \right|
\]
for $n \in \ZZ_{ >0 }$ extends to a polynomial in $n$ of degree $\dim(\cP)$.
(See, e.g.,~\cite{ccd} for background on Ehrhart theory.)
We will assume throughout that $\cP$ is full dimensional, and so the degree of $\ehr_\cP(n)$ is $d$.
An equivalent formulation of Ehrhart's theorem is that the \Def{Ehrhart series}
$
  1 + \sum_{ n \ge 1 } \ehr_\cP(n) \, z^n
$
evaluates to a rational function of the form $\frac{ h_\cP^*(z) }{ (1-z)^{ d+1 } }$ for some
polynomial $h_\cP^*(z)$ of degree $s \le d$, the \Def{$h^*$-polynomial} of $\cP$---a name
for the binomial transform of an Ehrhart polynomial that has become somewhat of a standard. 
We are interested in (linear) constraints among the $h^*$-coefficients.
Stanley~\cite{stanleydecomp} proved that the coefficients of $h_\cP^*(z)$ are
nonnegative integers. We will use below that the coefficient of $z^d$ equal the
number of interior lattice points of $\cP$, and the constant terms equals~1.

The \Def{Ehrhart--Macdonald reciprocity theorem}~\cite{macdonald} gives the algebraic relation
\[
  (-1)^d \ehr_\cP(-n) \ = \ \ehr_{ \cP^\circ } (n)
\]
where $\cP^\circ$ denotes the interior of $\cP$. An equivalent version is
\begin{equation}\label{eq:ehrmacrec}
  z^{ d+1 } h^*_\cP( \tfrac 1 z ) \ = \ h^*_{ \cP^\circ } (z)
\end{equation}
where the \Def{$h^*$-polynomial} of $\cP^\circ$ is defined through
\[
  \sum_{ n \ge 1 } \ehr_{ \cP^\circ }(n) \, z^n = \frac{ h^*_{ \cP^\circ } (z) }{
(1-z)^{ d+1 } } \, .
\]
Note that the degree of $h^*_{ \cP^\circ } (z)$ equals $d+1$.

A \emph{triangulation} of a $d$-dimensional polytope $\cP$ is a collection of simplices so that their union is $\cP$ and the intersection of
two simplices is a face of both. 
(See, e.g.,~\cite{deloerarambausantos} for background on triangulations.)
A triangulation of $\cP$ is \emph{unimodular} if all simplices have integer vertices and (minimal) volume $\frac 1 {d!}$.
A triangulation $T$ comes with an \Def{$f$-polynomial} 
\[
  f_T(z) \ := \ \sum_{j=0}^{d+1} f_{j-1} \, z^j
\]
where $f_j$ counts the number of $j$-dimensional faces of $T$ (and we set $f_{-1}=1$ for the empty face).
We further define the \Def{$h$-polynomial} of $T$ to be 
\[
h_T(z) \ := \ (1-z)^{d+1} \, f_T\left(\frac{z}{1-z}\right).
\]
If $\cP$ has a unimodular triangulation $T$, it is well known (see, e.g., \cite[Chapter~10]{ccd}) that the $h^*$-polynomial of $\cP$ equals the $h$-polynomial of~$T$.

Athanasiadis \cite[Theorem~1.3]{athanasiadishstareulerian} proved that, if $\cP$ is a $d$-dimensional
lattice polytope that admits a regular unimodular triangulation, then
\begin{align}
  &h^*_d \le h^*_{ d-1 } \le \dots \le h^*_{ \lfloor \frac{ d+1 }{ 2 } \rfloor } \label{eq:athan1} \\
  &h^*_{ j+1 } \ge h^*_{ d-j } \ \text{ for } \ 0 \le j \le \tfrac d 2 - 1 \label{eq:athan2} \\
  &h^*_j \le \binom{ h^*_1 + j - 1 }{ j } \ \text{ for } \ 0 \le j \le d \, . \label{eq:athan3}
\end{align}
Athanasiadis remarked in \cite{athanasiadishstareulerian} that these inequalities had been
independently proved by Hibi and Stanley (unpublished).
Stapledon~\cite[Theorem~2.20]{stapledondelta} showed that \eqref{eq:athan2} holds under the
(weaker) assumption that the boundary of $\cP$ admits a regular unimodular triangulation. 
Under the same condition, Stapledon proved that
\begin{equation}\label{eq:athan4}
  h^*_0 + \dots + h^*_{j+1} \le h^*_d + \dots + h^*_{ d-j } + \binom{ h^*_1 - h^*_d + j + 1 }{ j+1 } \ \text{ for } \ 0 \le j \le \tfrac d 2 - 1 \, .
\end{equation}
We remark that Stapledon derived \eqref{eq:athan2} and \eqref{eq:athan4} from a broad set of
$h^*$-inequalities, extending previous work of Betke--McMullen~\cite{betkemcmullen} and Payne~\cite{payneehrharttriang}.


\section{Order and Chromatic Polynomials}\label{sec:order}

Given a finite poset $(\Pi, \preceq)$ with $|\Pi| = d$, the \Def{order polynomial} $\Omega_\Pi^\circ(n)$ counts all strictly order-preserving maps from $\Pi$ to $[n] := \{ 1, 2, \dots, n \}$, i.e.,
\[
  \Omega_\Pi^\circ(n) \ := \ \left| \left\{ \varphi \in [n]^\Pi : \, a \prec b \ \Longrightarrow \varphi(a) < \varphi(b) \right\} \right| .
\]
Order polynomials first surfaced in~\cite{stanleychromaticlike}; we will encode them via
\[
  \sum_{ n \ge 1 } \Omega_\Pi^\circ(n) \, z^n \ = \ \frac{ \Omega^*_\Pi(z) }{ (1-z)^{ d+1 } } \, .
\]
(See, e.g.,~\cite{stanleyec1} for background on posets and order polynomials.)
Order polynomials are Ehrhart polynomials in disguise. We define the \Def{order polytope} of $\Pi$ as
\[
  \cO \ := \ \left\{ \varphi \in [0,1]^\Pi : \, a \preceq b \ \Longrightarrow \ \varphi(a) \le \varphi(b) \right\} .
\]
This much-studied subpolytope of the unit cube in $\RR^\Pi$ was introduced in~\cite{stanleyposetpolytopes}. From its definition we deduce that
\[
  \Omega_\Pi^\circ(n) \ = \ \ehr_{ \cO^\circ } (n+1) \, .
\]
This implies $\Omega^*_\Pi(z) = \frac 1 z \, h^*_{ \cO^\circ } (z) = z^d \, h^*_\cO (\frac
1 z)$, i.e.,
\begin{equation}\label{eq:htoomega}
  \Omega^*_j \ = \ h^*_{ d-j } 
\end{equation}
where the numbers on the right-hand side are the coefficients of $h^*_\cO(z)$.
Note that $\Omega^*_0 = h^*_d = 0$ (because $\cO$ contains no interior lattice points) and
$\Omega^*_d = h^*_0 = 1$.

\begin{theorem}\label{thm:orderpoldecomp}
Let $\Pi$ be a poset on $d$ elements and, as above, denote the binomial transform of its
order polynomial by $\Omega^*_\Pi(z) = \Omega^*_d \, z^d + \Omega^*_{ d-1 } \, z^{ d-1 } +
\dots + \Omega^*_1 \, z$.
Then
\begin{align*}
  &\Omega^*_{ 1 } \le \Omega^*_{ 2 } \le \dots \le \Omega^*_{ \lfloor \frac{ d+1 }{ 2 } \rfloor } \\
  &\Omega^*_{ j } \le \Omega^*_{ d-j } \ \text{ for } \ 1 \le j \le \tfrac {d-1} 2 \\
  &\Omega^*_{d-j} \le \binom{ \Omega^*_{d-1} + j - 1 }{ j } \ \text{ for } \ 0 \le j \le d-1 \\
  &\Omega^*_d + \dots + \Omega^*_{d-j} \le \Omega^*_{1} + \dots + \Omega^*_{ j } + \binom{
\Omega^*_{ d-1 } - \Omega^*_{1} + j }{ j } \ \text{ for } \ 1 \le j \le \tfrac {d-1} 2 \, . 
\end{align*}
\end{theorem}

\begin{proof} 
Let $\mu: \RR^d \to H_0 := \left\{ \bx \in \RR^d : \, x_1 + x_2 + \dots + x_d = 0 \right\}$ be an orthogonal projection, and let
$L$ be the lattice in $H_0$ generated by $\mu(\be_1), \mu(\be_2), \dots, \mu(\be_{d-1})$, i.e., $L = \mu(\ZZ^d)$, and $\mu(\be_d) = -
\mu(\be_1) - \dots - \mu(\be_{ d-1 })$.

We claim that the order polytope $\cO$ of $\Pi$ and $\mu(\cO)$ have the same $h^*$-polynomial.
To see this, consider the canonical unimodular triangulation $T$ of $\cO$, using the hyperplanes $x_j = x_k$.
The image of each simplex $\Delta \in T$ under the projection $\mu$ is a unimodular simplex in $H_0$ (with respect to $L$), and the vertices
$(0, \dots, 0)$ and $(1, \dots, 1)$ both get projected to the origin. This gives a unimodular triangulation $T_\mu$ of $\mu(\cO)$, and $T$ is
combinatorially a cone over $T_\mu$; in particular, the $f$-vectors of $T$ and $T_\mu$ are related via
\[
  f_T(z) \ = \ f_{ T_\mu } (z) \, (1+z) \, .
\]
Because both triangulations are unimodular,\footnote{
The equality~\eqref{eq:sameh*} of $h^*_\cO (z)$ and $h^*_{ \mu(\cO) } (z)$ can be also seen by noticing that the triangulations $T$ and
$T_\mu$ are regular and therefore shellable, and they have the same $h$-polynomial. See, e.g., \cite{haaseetalunimtriangpos} why order polytopes are compressed, and therefore have regular unimodular triangulations, and also how these properties are preserved under the projection $\mu$.
We also note that projected order polytopes are examples of \emph{alcoved polytopes}~\cite{lampostnikov}.
}
\begin{equation}\label{eq:sameh*}
\begin{aligned}
  h^*_\cO (z)
  \ &= \ h_T(z)
  \  = \ (1-z)^{ d+1 } f_T \left( \frac{ z }{ 1-z } \right) 
  \  = \ (1-z)^{ d+1 } \left( 1 + \frac{ z }{ 1-z } \right) f_{T_\mu} \left( \frac{ z }{ 1-z } \right) \\
  \ &= \ (1-z)^{ d } \, f_{T_\mu} \left( \frac{ z }{ 1-z } \right) 
  \  = \ h_{ T_\mu } (z)
  \  = \ h^*_{ \mu(\cO) } (z) \ . 
\end{aligned}
\end{equation}
The coefficients of $h^*_{ \mu(\cO) } (z)$ satisfy the Athanasiadis--Stapledon inequalities
\eqref{eq:athan1}--\eqref{eq:athan4}, with $d$ replaced by $d-1$.
Via \eqref{eq:htoomega}, this implies
\begin{align*}
  &\Omega^*_{ 1 } \le \Omega^*_{ 2 } \le \dots \le \Omega^*_{ \lfloor \frac{ d+1 }{ 2 } \rfloor } \\
  &\Omega^*_{ j } \le \Omega^*_{ d-j } \ \text{ for } \ 1 \le j \le \tfrac {d-1} 2 \\
  &\Omega^*_{d-j} \le \binom{ \Omega^*_{d-1} + j - 1 }{ j } \ \text{ for } \ 0 \le j \le d-1 \\
  &\Omega^*_d + \dots + \Omega^*_{d-j} \le \Omega^*_{1} + \dots + \Omega^*_{ j } + \binom{
\Omega^*_{ d-1 } - \Omega^*_{1} + j }{ j } \ \text{ for } \ 1 \le j \le \tfrac {d-1} 2 \, . \qedhere
\end{align*}
\end{proof}

\begin{proof}[Proof of the first 
two sets of inequalities in Theorem~\ref{thm:chrompoldecomp}]
Let $A(G)$ be the set of all acyclic orientations of $G$.\footnote{
An orientation is \Def{acyclic} if it does not contain any coherently directed cycles.
}
Then the chromatic polynomial $\chi_G(n)$ of $G$ decomposes naturally into order polynomials as
\begin{equation}\label{eq:chiintoorder}
  \chi_G(n) \ = \sum_{ \Pi \in A(G) } \Omega_\Pi^\circ (n) \, .
\end{equation}
Here we identify an acyclic orientation $\Pi$ with its corresponding poset.
(In this language, it is quite natural to think of $\Omega_\Pi^\circ (n)$ as the chromatic polynomial of the digraph $\Pi$.)
Because every $\Pi \in A(G)$ has $d$ elements,
\begin{equation*} 
  \chi^*_G (z) \ = \ \sum_{ \Pi \in A(G) } \Omega^*_\Pi(z) \, ,
\end{equation*}
and so the first 
two sets of inequalities in Theorem~\ref{thm:chrompoldecomp} follow from Theorem~\ref{thm:orderpoldecomp}.
\end{proof}


\section{Flow and tension polynomials}\label{sec:flow}

The inequalities for the coefficients of $\phi^*_G (z)$, $f^*_G (z)$, and $t^*_G (z)$ are proved in a
similar way, except that now we do not have the luxury of the dimension reduction exhibited
in the proof of Theorem~\ref{thm:orderpoldecomp}.

\begin{proof}[Proof of the remaining inequalities in Theorem~\ref{thm:chrompoldecomp}]
We start by showing that each of $\phi^*_G (z)$, $f^*_G (z)$, and $t^*_G (z)$ is the sum of
$h^*$-polynomials of open polytopes of the same dimension. (By an \emph{open polytope} we
simply mean the interior of a polytope.)

\begin{itemize}

\item For $\ZZ_n$-flows we use \cite[Proposition~2.3]{breuersanyal}, which expresses $\phi_G(n)$ as a sum of Ehrhart polynomials of certain open polytopes, all of which have dimension $\xi$.
(Briefly, one replaces the flow equations over $\ZZ_n$ by a set of affine equations over
$\RR$, in which $n$ now acts as a dilation parameter.)
Thus $\phi^*_G(z)$ is a sum of $h^*$-polynomials of open polytopes of dimension~$\xi$.

\item For integer flows, we write, as in the proof of \cite[Theorem~1]{kocholpolynomial},
\[
  f_G(n) \ = \sum_{ \Pi \in T(G) } p_\Pi(n)
\]
where $T(G)$ is the set of all totally cyclic orientations of $G$,\footnote{
An orientation is \Def{totally cyclic} if every edge lies in a coherently directed cycle.
} 
and $p_\Pi(n)$ counts the $\ZZ$-flows $x$ on $\Pi$ whose images satisfy $0< x(e) < n$. As noted in \cite{kocholpolynomial}, $p_\Pi(n)$ is the Ehrhart polynomial of an open polytope with dimension $\xi$, and so
$f^*_G(z)$ is a sum of $h^*$-polynomials of open polytopes.

\item Similarly, for integer tensions, we use \cite[Section 4]{kocholtension} to write
\[
  t_G(n) \ = \sum_{ \Pi \in A(G) } u_\Pi(n)
\]
where, as above, $A(G)$ is the set of all acyclic orientations of $G$, and $u_\Pi(n)$ counts the
$\ZZ$-tensions $x$ on $\Pi$ with $0< x(e) < n$. By \cite{kocholtension}, $u_\Pi(n)$ is the
Ehrhart polynomial of an open polytope with dimension $d-c$, and so $t^*_G(z)$ is a sum of $h^*$-polynomials of open polytopes.

\end{itemize}

In each of the above three cases, the polytopes in the decomposition admit regular
unimodular triangulations (see, e.g., \cite{haaseetalunimtriangpos,breuerdall}), so we can
apply \eqref{eq:athan1}--\eqref{eq:athan4} and these inequalities will then extend linearly to 
the coefficients of $\phi^*_G (z)$, $f^*_G (z)$, and $t^*_G (z)$.

It remains to rewrite \eqref{eq:athan1}--\eqref{eq:athan4} for $h^*_{ \cP^\circ } (z)
= \alpha_{ d+1 } \, z^{ d+1 } + \alpha_{ d } \, z^{ d } + \dots + \alpha_1 \, z$ (assuming
the polytope in question has dimension $d$) via \eqref{eq:ehrmacrec}:
\begin{align*}
  &\alpha_1 \le \alpha_{ 2 } \le \dots \le \alpha_{ \lfloor \frac{ d }{ 2 } \rfloor + 1 } \\
  &\alpha_{ d-j } \ge \alpha_{ j+1 } \ \text{ for } \ 0 \le j \le \tfrac d 2 - 1 \\
  &\alpha_{d+1-j} \le \binom{ \alpha_d + j - 1 }{ j } \ \text{ for } \ 0 \le j \le d \\
  &\alpha_{d+1} + \dots + \alpha_{d-j} \le \alpha_1 + \dots + \alpha_{ j+1 } + \binom{
\alpha_d - \alpha_1 + j + 1 }{ j+1 } \ \text{ for } \ 0 \le j \le \tfrac d 2 - 1 \, .
\end{align*}
The remaining inequalities in Theorem~\ref{thm:chrompoldecomp} now follow from the first two
sets of inequalities above.
\end{proof}


\bibliographystyle{amsplain}

\begin{thebibliography}{10}

\bibitem{athanasiadishstareulerian}
Christos~A. Athanasiadis, \emph{{$h^\ast$}-vectors, {E}ulerian polynomials and
  stable polytopes of graphs}, Electron. J. Combin. \textbf{11} (2004/06),
  no.~2, Research Paper 6, 13 pp.

\bibitem{ccd}
Matthias Beck and Sinai Robins, \emph{Computing the {C}ontinuous {D}iscretely:
  Integer-point {E}numeration in {P}olyhedra}, second ed., Undergraduate Texts
  in Mathematics, Springer, New York, 2015, electronically available at {\tt
  http://math.sfsu.edu/beck/ccd.html}.

\bibitem{betkemcmullen}
Ulrich Betke and Peter McMullen, \emph{Lattice points in lattice polytopes},
  Monatsh. Math. \textbf{99} (1985), no.~4, 253--265.

\bibitem{brentichromatic}
Francesco Brenti, \emph{Expansions of chromatic polynomials and log-concavity},
  Trans. Amer. Math. Soc. \textbf{332} (1992), no.~2, 729--756.

\bibitem{breuerdall}
Felix Breuer and Aaron Dall, \emph{Bounds on the coefficients of tension and
  flow polynomials}, J. Algebraic Combin. \textbf{33} (2011), no.~3, 465--482,
  {\tt arXiv:1004.3470}.

\bibitem{breuersanyal}
Felix Breuer and Raman Sanyal, \emph{Ehrhart theory, modular flow reciprocity,
  and the {T}utte polynomial}, Math. Z. \textbf{270} (2012), no.~1-2, 1--18,
  {\tt arXiv:0907.0845}.

\bibitem{deloerarambausantos}
Jes{\'u}s~A. De~Loera, J{\"o}rg Rambau, and Francisco Santos,
  \emph{Triangulations}, Algorithms and Computation in Mathematics, vol.~25,
  Springer-Verlag, Berlin, 2010.

\bibitem{ehrhartpolynomial}
Eug{\`e}ne Ehrhart, \emph{Sur les poly\`edres rationnels homoth\'etiques \`a
  {$n$}\ dimensions}, C. R. Acad. Sci. Paris \textbf{254} (1962), 616--618.

\bibitem{gansnervo}
Emden~R. Gansner and Kiem~Phong Vo, \emph{The chromatic generating function},
  Linear and Multilinear Algebra \textbf{22} (1987), no.~1, 87--93.

\bibitem{gesselchromatic}
Ira~M. Gessel, \emph{Acyclic orientations and chromatic generating functions},
  Discrete Math. \textbf{232} (2001), no.~1-3, 119--130.

\bibitem{haaseetalunimtriangpos}
Christian Haase, Andreas Paffenholz, Lindsay~C. Piechnik, and Francisco Santos,
  \emph{Existence of unimodular triangulations---positive results}, Preprint
  ({\tt arXiv:1405.1687}).

\bibitem{hershswartz}
Patricia Hersh and Ed~Swartz, \emph{Coloring complexes and arrangements}, J.
  Algebraic Combin. \textbf{27} (2008), no.~2, 205--214, {\tt
  arXiv:math/0706.3657}.

\bibitem{jaegernwzflowproblems}
Fran\c{c}ois Jaeger, \emph{Nowhere-zero flow problems}, Lowell W. Beineke,
  Robin J. Wilson (Eds.), Selected topics in graph theory, vol. 3, Academic
  Press, San Diego, CA, 1988, pp.~71--95.

\bibitem{kocholpolynomial}
Martin Kochol, \emph{Polynomials associated with nowhere-zero flows}, J.
  Combin. Theory Ser. B \textbf{84} (2002), no.~2, 260--269.

\bibitem{kocholtension}
\bysame, \emph{Tension polynomials of graphs}, J. Graph Theory \textbf{40}
  (2002), no.~3, 137--146.

\bibitem{korfhage}
Robert~R. Korfhage, \emph{{$\sigma $}-polynomials and graph coloring}, J.
  Combinatorial Theory Ser. B \textbf{24} (1978), no.~2, 137--153.

\bibitem{lampostnikov}
Thomas Lam and Alexander Postnikov, \emph{Alcoved polytopes. {I}}, Discrete
  Comput. Geom. \textbf{38} (2007), no.~3, 453--478, {\tt arXiv:math/0501246}.

\bibitem{leon}
Emerson Le\'{o}n, \emph{Stapledon decompositions and inequalities for
  coefficients of chromatic polynomials}, S\'{e}m. Lothar. Combin. \textbf{78B}
  (2017), Art. 24, 12.

\bibitem{linialchromatic}
Nathan Linial, \emph{Graph coloring and monotone functions on posets}, Discrete
  Math. \textbf{58} (1986), no.~1, 97--98.

\bibitem{macdonald}
Ian~G. Macdonald, \emph{Polynomials associated with finite cell-complexes}, J.
  London Math. Soc. (2) \textbf{4} (1971), 181--192.

\bibitem{payneehrharttriang}
Sam Payne, \emph{Ehrhart series and lattice triangulations}, Discrete Comput.
  Geom. \textbf{40} (2008), no.~3, 365--376, {\tt arXiv:math/0702052}.

\bibitem{seymourhandbook}
Paul~D. Seymour, \emph{Nowhere-zero flows}, Handbook of combinatorics, {V}ol.\
  1,\ 2, Elsevier Sci. B. V., Amsterdam, 1995, Appendix: Colouring, stable sets
  and perfect graphs, pp.~289--299.

\bibitem{stanleychromaticlike}
Richard~P. Stanley, \emph{A chromatic-like polynomial for ordered sets}, Proc.
  {S}econd {C}hapel {H}ill {C}onf. on {C}ombinatorial {M}athematics and its
  {A}pplications ({U}niv. {N}orth {C}arolina, {C}hapel {H}ill, {N}.{C}., 1970),
  Univ. North Carolina, Chapel Hill, N.C., 1970, pp.~421--427.

\bibitem{stanleydecomp}
\bysame, \emph{Decompositions of rational convex polytopes}, Ann. Discrete
  Math. \textbf{6} (1980), 333--342.

\bibitem{stanleyposetpolytopes}
\bysame, \emph{Two poset polytopes}, Discrete Comput. Geom. \textbf{1} (1986),
  no.~1, 9--23.

\bibitem{stanleyec1}
\bysame, \emph{Enumerative {C}ombinatorics. {V}olume 1}, second ed., Cambridge
  Studies in Advanced Mathematics, vol.~49, Cambridge University Press,
  Cambridge, 2012.

\bibitem{stapledondelta}
Alan Stapledon, \emph{Inequalities and {E}hrhart {$\delta$}-vectors}, Trans.
  Amer. Math. Soc. \textbf{361} (2009), no.~10, 5615--5626, {\tt
  arXiv:math/0801.0873}.

\bibitem{tomescu}
Ioan Tomescu, \emph{Graphical {E}ulerian numbers and chromatic generating
  functions}, Discrete Math. \textbf{66} (1987), no.~3, 315--318.

\bibitem{tutteflowpoly}
William~T. Tutte, \emph{A ring in graph theory}, Proc. Cambridge Philos. Soc.
  \textbf{43} (1947), 26--40.

\bibitem{wilfchromatic}
Herbert~S. Wilf, \emph{Which polynomials are chromatic?}, Colloquio
  {I}nternazionale sulle {T}eorie {C}ombinatorie ({R}oma, 1973), {T}omo {I},
  Accad. Naz. Lincei, Rome, 1976, pp.~247--256.

\end{thebibliography}

\def\cprime{$'$} \def\cprime{$'$}
\providecommand{\bysame}{\leavevmode\hbox to3em{\hrulefill}\thinspace}
\providecommand{\MR}{\relax\ifhmode\unskip\space\fi MR }
\providecommand{\MRhref}[2]{%
  \href{http://www.ams.org/mathscinet-getitem?mr=#1}{#2}
}
\providecommand{\href}[2]{#2}

\setlength{\parskip}{0cm} 

\end{document}